\documentclass[aap,secthm,MSNbibl,citesort,dvips]{arximspdf}

\doi{10.1214/10-AAP686}
\volume{20}
\issue{6}
\pubyear{2010}
\firstpage{2295}
\lastpage{2317}

\makeatletter

\def\sfrac#1#2{#1/#2}

\def\afrac#1#2{#1/(#2)}

\newtheorem{theorem}{Theorem}[section]
\newtheorem{lemma}{Lemma}[section]
\newtheorem{corollary}{Corollary}[section]
\newtheorem{proposition}{Proposition}[section]

\newproclaim{definition}{Definition}[section]
\newproclaim{remark}{Remark}[section]
\newproclaim{convention}{Convention}[section]
\newproclaim{claim}{Claim}[section]
\newproclaim{case}{Case}
\newproclaim{casee}{Case}
\newproclaim{caseee}{Case}
\newproclaim{subcase}{Subcase}

\makeatother

\begin{document}
\begin{frontmatter}

\title{Stationary systems of Gaussian processes}
\runtitle{Stationary systems of Gaussian processes}

\begin{aug}
\author{\fnms{Zakhar} \snm{Kabluchko}\corref{}\ead[label=e1]{kabluch@math.uni-goettingen.de}}
\runauthor{Z. Kabluchko}
\affiliation{Georg-August-Universit\"at G\"ottingen}

\address{Institut f\"ur Mathematische Stochastik\\
Georg-August-Universit\"at G\"ottingen\\
Goldschmidtstr. 7\\
D-37077 G\"ottingen\\
Germany\\
\printead{e1}}
\end{aug}

\received{\smonth{3} \syear{2009}}
\revised{\smonth{12} \syear{2009}}

%
\begin{abstract}
We describe all countable particle systems on $\mathbb{R}$ which have the
following three properties:
independence, Gaussianity and stationarity. More precisely,
we consider particles on the real line starting at the points of a
Poisson point process with
intensity measure $\mathfrak{m}$ and moving independently of each other
according to the law of
some Gaussian process $\xi$. We classify all pairs $(\mathfrak{m},\xi)$
generating a \textit{stationary}
particle system, obtaining three families of examples. In the first,
trivial family, the measure
$\mathfrak{m}$ is arbitrary, whereas the process $\xi$ is stationary.
In the
second family, the measure
$\mathfrak{m}$ is a multiple of the Lebesgue measure, and $\xi$ is
essentially a
Gaussian stationary
increment process with linear drift. In the third, most interesting
family, the measure
$\mathfrak{m}$ has a density of the form $\alpha e^{-\lambda x}$,
where $\alpha
>0$, $\lambda\in\mathbb{R}$,
whereas the process $\xi$ is of the form $\xi(t)=W(t)-\lambda\sigma
^2(t)/2+c$, where $W$
is a zero-mean Gaussian process with stationary increments, $\sigma
^2(t)=\operatorname{Var}W(t)$, and $c\in\mathbb{R}$.
\end{abstract}

%

\end{frontmatter}

\section{Introduction}\label{sec:intro}

\subsection{Statement of the problem}\label{subsec:problem}

Stationary systems of particles evolving independently of each other
according to the law of a Markov process have been extensively studied
by many authors (see, e.g., the monographs~\cite{masipresuttibook},
Chapter~1, \cite{kipnislandimbook},
Chapter~1, \cite{reveszbook},
as well as the papers~\cite
{derman55,brown70,martinloef76,liggettport88,coxgriffeath84,deuschelwang94},
to cite only a few references).
The aim of the present paper is to study systems of particles evolving
independently of each other in a \textit{Gaussian} rather than \textit
{Markovian} way. Our main result provides a classification of all those
Gaussian particle systems which are \textit{stationary}.

We are interested in at most countable systems of particles moving
randomly on the real line in such a way that the following three
requirements are satisfied:
\begin{enumerate}[(A1)]
\item[(A1)] The particles are independent of each other.
\item[(A2)] The law describing the motion of each particle is Gaussian
and the same for all particles.
\item[(A3)] The particles are in an equilibrium.
\end{enumerate}

The independence stated in requirement~(A1) implies that the starting
positions of particles should be scattered independently over $\mathbb{R}$,
which, in more rigorous terms, means that they should form a not
necessarily homogeneous Poisson point process on $\mathbb{R}$. Requirement~(A2)
means that the stochastic processes describing the deviations of the
particles from their starting positions should be Gaussian, having the
same law for all particles, and, by requirement~(A1), independent of
each other.

In view of this, the meaning of the first two requirements may be
described in rigorous terms as follows.
Let $\{U_i, i\in\mathbb{N}\}$ be a Poisson point process on $\mathbb
{R}$ with intensity
measure $\mathfrak{m}$. We will always assume that $\mathfrak{m}$
satisfies the following
integrability condition:
%
\begin{equation}\label{eq:mm_cond_fin}
\int_{\mathbb{R}}e^{-\varepsilon x^2} \mathfrak{m}(dx)<\infty\qquad\mbox{for every } \varepsilon>0.
\end{equation}
In most cases of interest, the measure $\mathfrak{m}$ will be
infinite, and so
let us agree to use $\mathbb{N}$ as an index set for the points $U_i$, even
though the case where $\mathfrak{m}$ is finite (and, hence, a.s.\ only finitely
many points $U_i$ exist) is not formally excluded.

Let $\xi_i, i\in\mathbb{N}$, be independent copies of a Gaussian
process $\{
\xi
(t), t\in\mathbb{R}^d\}$.
We define $V_i(t)$, the position of $i$th particle at time $t\in
\mathbb{R}^d$
(which we allow to be multidimensional), by
%
\begin{equation}\label{eq:def_v}
V_i(t)=U_i+\xi_i(t).
\end{equation}

\begin{definition}
The random collection of functions $\mathfrak{P}=\{V_i, i\in\mathbb
{N}\}$ will be
called the \textit{independent Gaussian particle system} (or simply
\textit{Gaussian system}) generated by the pair $(\mathfrak{m},\xi
)$. We use the
notation $\mathit{GS}(\mathfrak{m},\xi)$.
\end{definition}
%
\begin{remark}
It should be stressed that we \textit{do not} assume the process $\xi$
to have zero mean, which means that we allow for a deterministic
component in the random motion of particles.
In general, it also may happen that $\xi(0)\neq0$, in which case the
particles make nonzero jumps immediately after starting at $U_i$.
\end{remark}

Let us turn to requirement~(A3). Given $t_1,\ldots,t_n\in\mathbb
{R}^d$, we
define a point process $\mathfrak{P}_{t_1,\ldots,t_n}$ on $\mathbb
{R}^n$ by recording
the positions of particles at times $t_1,\ldots,t_n$. That is, we set
%
\begin{equation}\label{eq:def_ppp_fin_dim}
\mathfrak{P}_{t_1,\ldots,t_n}=\{(V_i(t_1),\ldots,V_i(t_n)), i\in
\mathbb{N}\}.
\end{equation}
The family $\{\mathfrak{P}_{t_1,\ldots,t_n}\dvtx  n\in\mathbb{N},
t_1,\ldots,t_n\in\mathbb{R}
^d\}$
may be viewed as the family of ``finite-dimensional distributions'' of
$\mathfrak{P}$.
%
\begin{definition}\label{def:gauss_stat}
A Gaussian system $\mathfrak{P}$ is called \textit{stationary}
if for every $n\in\mathbb{N}$, every $t_1,\ldots,t_n\in\mathbb
{R}^d$, and every
$h\in\mathbb{R}
^{d}$, we have the following equality of laws of point processes on
$\mathbb{R}^n$:
%
\begin{equation}\label{eq:stat_def}
\mathfrak{P}_{t_1+h,\ldots,t_n+h}\stackrel{d}{=}\mathfrak
{P}_{t_1,\ldots,t_n}.
\end{equation}
\end{definition}

The purpose of this paper is to provide a description of all stationary
Gaussian systems. Let us stress that for Markovian particle systems,
the corresponding question has a rather simple solution. Let the
initial positions of the particles be chosen to form a Poisson point
process with $\sigma$-finite intensity measure $\mathfrak{m}$ on some
measurable
space $(\Omega,\mathcal{A})$, and let the particles move
independently of each
other according to the law of some Markov process on $\Omega$ with
transition kernel $P(x,dy)$. Then by a result of~\cite{brown70}, the
particle system is stationary if and only if the measure $\mathfrak
{m}$ is
$P$-invariant (see also~\cite{doobbook}, page~404, and~\cite{derman55}, Theorem~2, for weaker results).

\subsection{Statement of the main result}
First we introduce some notation.
For $\lambda\in\mathbb{R}$, we denote by $\mathfrak{e}_{\lambda}$
a measure on $\mathbb{R}
$ with
a density of the form $e^{-\lambda x}$ with respect to the Lebesgue
measure. That is,
%
\begin{equation}
\mathfrak{e}_{\lambda}(dx)=e^{-\lambda x}\,dx.
\end{equation}
In particular, $\mathfrak{e}_0$ is the Lebesgue measure itself.

A function $f\dvtx \mathbb{R}^d\to\mathbb{R}$ is called \textit
{additive} if
$f(t_1+t_2)=f(t_1)+f(t_2)$ for every $t_1,t_2\in\mathbb{R}^d$. Under minor
additional assumptions, say, measurability, an additive function must
be of the form $f(t)=\langle c,t\rangle$ for some $c\in\mathbb{R}^d$.
%
%
\begin{convention}
All stationary processes and processes with stationary increments are
always supposed to have \textit{zero mean}.
\end{convention}

The next theorem is our main result.
%
\begin{theorem}\label{theo:main_gen}
Let $\mathcal{S}$ be the set of all pairs $(\mathfrak{m},\xi)$,
where $\mathfrak{m}$ is a measure
satisfying~(\ref{eq:mm_cond_fin}) and $\{\xi(t), t\in\mathbb{R}^d\}
$ is a
Gaussian process, with the property that the particle system
$\mathit{GS}(\mathfrak{m}
,\xi
)$ is stationary. Then
%
\begin{equation}\label{eq:HH_decomp}
\mathcal{S}=\mathcal{S}_1\cup\mathcal{S}_2\cup\mathcal{S}_3,
\end{equation}
where the sets $\mathcal{S}_1,\mathcal{S}_2,\mathcal{S}_3$ are
defined as follows:
\begin{enumerate}
\item The set $\mathcal{S}_1$ consists of all pairs $(\mathfrak
{m},\xi)$, where $\mathfrak{m}
$ is
an arbitrary measure on $\mathbb{R}$ satisfying~(\ref
{eq:mm_cond_fin}), and
\[
\{\xi(t), t\in\mathbb{R}^d\}\stackrel{d}{=}\{W(t)+c, t\in\mathbb
{R}^d\}
\]
for some stationary Gaussian process $\{W(t), t\in\mathbb{R}^d\}$ and some
$c\in
\mathbb{R}$.
\item The set $\mathcal{S}_2$ consists of all pairs $(\mathfrak
{m},\xi)$, where
\[
\mathfrak{m}=\alpha\mathfrak{e}_0 \quad\mbox{and}\quad \{\xi(t), t\in
\mathbb{R}^d\}\stackrel{d}{=}\{
W(t)+f(t)+c, t\in\mathbb{R}^d\}
\]
for some $\alpha>0$, $c\in\mathbb{R}$, a Gaussian process $\{W(t),
t\in\mathbb{R}
^d\}$
with stationary increments, and an additive function $f\dvtx \mathbb
{R}^d\to\mathbb{R}$.
\item The set $\mathcal{S}_3$ consists of all pairs $(\mathfrak
{m},\xi)$, where
\[
\mathfrak{m}=\alpha\mathfrak{e}_{\lambda} \quad\mbox{and}\quad \{\xi(t),
t\in\mathbb{R}^d\}
\stackrel{d}{=}\{
W(t)-\lambda\sigma^2(t)/2+c, t\in\mathbb{R}^d\}
\]
for some $\alpha>0$, $\lambda\neq0$, $c\in\mathbb{R}$, and some
Gaussian process
$\{W(t),t\in\mathbb{R}^d\}$ with stationary increments and variance
$\sigma^2(t)$.
\end{enumerate}
\end{theorem}




The stationarity of Gaussian systems of type $\mathcal{S}_1$ is a rather
trivial fact and is due to the stationarity of the driving process $\xi
$. Somewhat less trivial, but still rather appealing, is the fact that
Gaussian systems of type $\mathcal{S}_2$ are stationary. An example of a
Gaussian system of type $\mathcal{S}_2$ can be obtained by taking
$\mathfrak{m}$ to be
the Lebesgue measure on $\mathbb{R}$ and $\xi$ to be a (fractional) Brownian
motion with a linear drift.


Surprisingly, the class of stationary Gaussian systems is not exhausted
by the two ``trivial'' families $\mathcal{S}_1$ and $\mathcal{S}_2$:
there is one more,
nontrivial, family $\mathcal{S}_3$. An example of a Gaussian system
of type
$\mathcal{S}_3$ can be obtained by taking
\[
\mathfrak{m}=\mathfrak{e}_1 \quad\mbox{and}\quad \{\xi(t), t\in\mathbb
{R}\}\stackrel{d}{=}\{W_{\kappa
}(t)-|t|^{\kappa}, t\in\mathbb{R}\},
\]
where $\{W_{\kappa}(t), t\in\mathbb{R}\}$ is a fractional Brownian
motion with
index $\kappa\in(0,2]$, that is,\ a stationary increment Gaussian
process with
\[
\operatorname{Cov}(W_{\kappa}(t_1), W_{\kappa}(t_2))=|t_1|^{\kappa
}+|t_2|^{\kappa
}-|t_1-t_2|^{\kappa},\qquad t_1,t_2\in\mathbb{R}.
\]
For $\kappa=1$, this Gaussian system appeared in~\cite{brownresnick77}
in connection with maxima of independent Ornstein--Uhlenbeck processes.
For general $\kappa\in(0,2]$, the driving process $W_{\kappa
}(t)-|t|^{\kappa}$ appeared in~\cite{pickands69a}, also in connection
with maxima of Gaussian processes. In a similar way, particle systems
of type $\mathcal{S}_2$ appeared in~\cite{penrose91} in connection
with minima
(in the absolute value sense) of independent Gaussian processes. The
results of~\cite{brownresnick77} were generalized in~\cite
{kabluchkoetal}. In particular, it was shown in Theorem~2 of~\cite
{kabluchkoetal} that Gaussian systems of type $\mathcal{S}_3$ with an
additional requirement $\alpha=1$, $\lambda=1$, $c=0$ were stationary.
Gaussian systems of type $\mathcal{S}_3$ have some vague similarity
with the
``competing particle systems'' studied in~\cite{ruzmaikinaaizenman05}
(see also~\cite{arguinaizenman09,shkolnikov09}). Note that in contrast
to our setting, the particles in~\cite{ruzmaikinaaizenman05} evolve by
increments which are \textit{independent} in time.

At a first sight, it may look that the family $\mathcal{S}_2$ can be included
into the family $\mathcal{S}_3$ by allowing the parameter $\lambda$
in the
definition of $\mathcal{S}_3$ to be $0$. However, this is not the
case: the
family $\mathcal{S}_2$ has an additional ``degree of freedom''
represented by
the additive function $f$.


In view of particle systems interpretation of Theorem~\ref{theo:main_gen},
of special interest are stationary Gaussian systems
driven by a process $\xi$ satisfying $\xi(0)=0$. In the next corollary
we provide a classification of such systems, excluding for convenience
the noninteresting case in which $\xi$ is a version of the zero process.
%
\begin{corollary}\label{cor:main_gen}
Let $\mathfrak{m}$ be a measure satisfying~(\ref{eq:mm_cond_fin}),
and let $\{
\xi
(t), t\in\mathbb{R}^d\}$ be a Gaussian process with $\xi(0)=0$.
Assume that for
some $t_0$, $\xi(t_0)$ is not a.s.\ $0$. Then the particle system
$\mathit{GS}(\mathfrak{m},\xi)$ is stationary iff
$
\mathfrak{m}=\alpha\mathfrak{e}_{\lambda}
$
for some $\alpha>0$ and $\lambda\in\mathbb{R}$, and
\[
\{\xi(t), t\in\mathbb{R}^d\}\stackrel{d}{=}
\cases{
\{W(t)+f(t), t\in\mathbb{R}^d\}, &\quad\mbox{if }$\lambda=0$,\vspace*{2pt}\cr
\{W(t)-\lambda\sigma^2(t)/2, t\in\mathbb{R}^d\}, &\quad\mbox{if }$\lambda\neq0$,
}
\]
for some Gaussian process $\{W(t), t\in\mathbb{R}^d\}$ with stationary
increments, variance $\sigma^2(t)$, $W(0)=0$ and, eventually, an
additive function $f:\mathbb{R}^d\to\mathbb{R}$.
\end{corollary}



\subsection{Organization of the paper}
Our main result, Theorem~\ref{theo:main_gen}, will be proved in
Section~\ref{sec:class}.
Although Theorem~\ref{theo:main_gen} classifies all pairs $(\mathfrak
{m},\xi)$
generating a stationary Gaussian system, it does not tell how to decide
whether two given pairs $(\mathfrak{m}',\xi')$, $(\mathfrak{m}'',\xi
'')$ generate equal
in law Gaussian systems or not. This gap will be filled in Section~\ref
{sec:class_complete}.




\section{Proof of the main result}\label{sec:class}
\subsection{Idea of the proof}
In this section we prove Theorem~\ref{theo:main_gen}. The ``easy'' part
of Theorem~\ref{theo:main_gen} stating that Gaussian systems generated
by the pairs $(\mathfrak{m},\xi)\in\mathcal{S}_1\cup\mathcal
{S}_2\cup\mathcal{S}_3$ are stationary
will be established in Proposition~\ref{prop:stat_H1}. The proof of the
converse statement is much more difficult. The first step will be done
in Proposition~\ref{prop:proof_main_2}, where it is shown that a pair
$(\mathfrak{m},\xi)$ generating a stationary Gaussian system must
belong to
$\mathcal{S}
_1\cup\mathcal{S}_2\cup\mathcal{S}_3$ provided that the measure
$\mathfrak{m}$ is a linear
combination of the Lebesgue measure $\mathfrak{e}_0$ and a measure of
the form
$\mathfrak{e}_{\lambda}$. Such linear combinations are well behaved under
convolutions with Gaussian measures, which makes it possible to do
explicit calculations with one- and two-dimensional distributions of
$\mathit{GS}(\mathfrak{m},\xi)$. The second step, carried out in Section~\ref
{subsec:ident_m}, is to show that this additional assumption on the
measure $\mathfrak{m}$ is satisfied for \textit{most} (but not all!) pairs
$(\mathfrak{m}
,\xi)$ generating a stationary Gaussian system. Essentially, this is
done by applying a result of Deny~\cite{deny59} and several related
lemmas collected in Section~\ref{subsec:lem_conv} to the
one-dimensional distributions of $\mathit{GS}(\mathfrak{m},\xi)$. The pairs
for which the
additional assumption on $\mathfrak{m}$ is not satisfied are shown to
belong to
the family $\mathcal{S}_1$.

\subsection{Notation}\label{subsec:notation}
We start by introducing the notation. We always assume that $\mathfrak
{m}$ is a
measure on $\mathbb{R}$ satisfying the integrability condition~(\ref
{eq:mm_cond_fin}), and that $\{\xi(t), t\in\mathbb{R}^d\}$ is a Gaussian
process. The law of the process $\xi$ is uniquely determined by its
mean and covariance for which we use the notation
%
\begin{equation}\label{eq:def_mu}
\mu(t)=\mathbb E\xi(t), \qquad r(t_1,t_2)=\operatorname{Cov}(\xi(t_1),\xi(t_2)).
\end{equation}
Further, we define the variance and the incremental variance of $\xi$ by
%
\begin{equation}\label{eq:def_sigma}
\sigma^2(t)=\operatorname{Var}\xi(t), \qquad \gamma(t_1,t_2)=\operatorname
{Var}[\xi(t_1)-\xi(t_2)].
\end{equation}
We will often use the identity
%
\begin{equation}\label{eq:r_sigma_gamma}
r(t_1,t_2)=\tfrac12 \bigl(\sigma^2(t_1)+\sigma^2(t_2)-\gamma(t_1,t_2) \bigr).
\end{equation}
Given $t_1,\ldots, t_n\in\mathbb{R}^d$, the law of the random vector
$(\xi
(t_1),\ldots,\xi(t_n))$ is denoted by $\mathfrak{n}_{t_1,\ldots,t_n}$.

Let $\mathcal{B}(\mathbb{R}^n)$ be the Borel $\sigma$-algebra of
$\mathbb{R}^n$. For a set
$B\subset\mathbb{R}^n$ and $x\in\mathbb{R}$, it will be convenient
to define
\[
B-x=B-(x,\ldots,x).
\]
So, $B-x$ is obtained by shifting the set $B$ ``diagonally'' in the
direction of the vector $(1,\ldots,1)$.

Define $\mathfrak{P}_{t_1,\ldots,t_n}$, the finite-dimensional
distributions of
$\mathfrak{P}$, as in~(\ref{eq:def_ppp_fin_dim}). The transformation
theory of
Poisson point processes (see, e.g.,\ Proposition~3.8 in~\cite
{resnickbook87}) tells that $\mathfrak{P}_{t_1,\ldots,t_n}$ is a
Poisson point
process on $\mathbb{R}^n$ with intensity measure $\mathfrak
{m}_{t_1,\ldots,t_n}$ that is
defined by
%
\begin{equation}\label{eq:mm_t_R}
\hspace*{20pt}\mathfrak{m}_{t_1,\ldots,t_n}(B)=\int_{\mathbb{R}}\mathbb{P}[(\xi
(t_1),\ldots,\xi
(t_n))\in B-x]
\mathfrak{m}(dx),\qquad B\in\mathcal{B}(\mathbb{R}^n).
\end{equation}
In particular, we will often use that $\mathfrak{m}_t=\mathfrak
{m}*\mathfrak{n}_t$ for every
$t\in
\mathbb{R}^d$, where $*$ denotes the convolution of measures. Note that
condition~(\ref{eq:mm_cond_fin}) ensures that $\mathfrak
{m}_{t_1,\ldots,t_n}(B)$
is finite for every bounded $B\in\mathcal{B}(\mathbb{R}^n)$.

We can restate Definition~\ref{def:gauss_stat} as follows: A Gaussian
system $\mathfrak{P}$ is stationary if for every $n\in\mathbb{N}$,
every $t_1,\ldots,
t_n, h\in\mathbb{R}^d$, and every $B\in\mathcal{B}(\mathbb{R}^n)$,
%
\begin{equation}\label{eq:mm_eq_mm}
\mathfrak{m}_{t_1,\ldots,t_n}(B)=\mathfrak{m}_{t_1+h,\ldots,t_n+h}(B).
\end{equation}

We denote the one-dimensional Gaussian measure with expectation $\mu_0$
and variance $\sigma_0^2$ by $\mathfrak{N}(\mu_0,\sigma_0^2)$. For future
reference, let us recall the following formula for the Laplace
transform of a Gaussian distribution:
%
\begin{equation}\label{eq:laplace_gauss}
\mbox{ if } N\sim\mathfrak{N}(\mu_0,\sigma_0^2),\qquad
\mbox{then } \mathbb Ee^{y N}=e^{\mu_0 y+ \sfrac{\sigma_0^2 y^2}{2}}.
\end{equation}

\subsection{\texorpdfstring{Proof of the easy part of Theorem~\protect\ref
{theo:main_gen}}{Proof of the easy part of Theorem~1.1}}
\label{subsec:easy_part}

In the next proposition we prove that Gaussian systems of types
$\mathcal{S}
_1,\mathcal{S}_2,\mathcal{S}_3$ are indeed stationary.
%
\begin{proposition}\label{prop:stat_H1}
Let $\mathfrak{P}=\mathit{GS}(\mathfrak{m},\xi)$, where $(\mathfrak{m},\xi
)\in\mathcal{S}_1\cup\mathcal{S}_2\cup
\mathcal{S}_3$.
Then $\mathfrak{P}$ is stationary.
\end{proposition}
\begin{pf}
Suppose that $(\mathfrak{m},\xi)\in\mathcal{S}_1$. By definition of
$\mathcal{S}_1$, we have
the following equality of laws, valid for all $n\in\mathbb{N}$,
$t_1,\ldots
,t_n,h\in\mathbb{R}^d$:
\[
(\xi(t_i))_{i=1}^n \stackrel{d}{=}\bigl(\xi(t_i+h)\bigr)_{i=1}^n.
\]
Let $B\subset\mathbb{R}^n$ be any Borel set. By~(\ref{eq:mm_t_R}),
we have
\begin{eqnarray*}
\mathfrak{m}_{t_1,\ldots,t_n}(B)
&=&\int_{\mathbb{R}}\mathbb{P}[(\xi(t_1),\ldots,\xi(t_n))\in
B-z]\mathfrak{m}(dz)\\
&=&\int_{\mathbb{R}}\mathbb{P}\bigl[\bigl(\xi(t_1+h),\ldots,\xi(t_n+h)\bigr)\in
B-z\bigr]\mathfrak{m}(dz)\\
&=&\mathfrak{m}_{t_1+h,\ldots,t_n+h}(B).
\end{eqnarray*}
Hence, equation~(\ref{eq:mm_eq_mm}) holds and $\mathfrak{P}$ is stationary.

Suppose that $(\mathfrak{m},\xi)\in\mathcal{S}_2$.
By definition of $\mathcal{S}_2$, we have $\mathfrak{m}=\alpha
\mathfrak{e}_0$ for some
$\alpha
>0$, and
%
\begin{equation}\label{eq:xi_stat_incr}
\bigl(\xi(t_i)-\xi(t_1)\bigr)_{i=1}^n \stackrel{d}{=}\bigl(\xi(t_i+h)-\xi(t_1+h)\bigr)_{i=1}^n
\end{equation}
for all $n\in\mathbb{N}$, $t_1,\ldots,t_n, h\in\mathbb{R}^d$. Let
$B\subset\mathbb{R}^n$ be
any Borel set.
Using~(\ref{eq:mm_t_R}) and~(\ref{eq:xi_stat_incr}), we obtain
%
\begin{eqnarray}\label{eq:wspom_br_br}
&&\mathfrak{m}_{t_1,\ldots,t_n}(B)\nonumber\\
&&\qquad=\alpha\int_{\mathbb{R}}\int_{\mathbb{R}^n} 1_{B-x}(y_1,\ldots,y_n) \mathfrak{n}_{t_1,\ldots,t_n}(dy_1,\ldots,dy_n)\,dx \nonumber\\
&&\qquad=\alpha\int_{\mathbb{R}}\int_{\mathbb{R}^n}1_{B-(x+y_1)}(0,y_2-y_1,\ldots,y_n-y_1)\mathfrak{n}_{t_1,\ldots,t_n}(dy_1,\ldots,dy_n)\,dx \nonumber\\
&&\qquad=\alpha\int_{\mathbb{R}}\int_{\mathbb{R}^n}1_{B-z}(0,y_2-y_1,\ldots,y_n-y_1)\mathfrak{n}_{t_1,\ldots,t_n}(dy_1,\ldots,dy_n)\, dz \\
&&\qquad=\alpha\int_{\mathbb{R}} \mathbb{P}\bigl[\bigl(\xi(t_i)-\xi(t_1)\bigr)_{i=1}^n\in B-z\bigr]\,dz\nonumber\\
&&\qquad=\alpha\int_{\mathbb{R}} \mathbb{P}\bigl[\bigl(\xi(t_i+h)-\xi(t_1+h)\bigr)_{i=1}^n \in B-z\bigr]\,dz\nonumber\\
&&\qquad=\mathfrak{m}_{t_1+h,\ldots,t_n+h}(B).\nonumber
\end{eqnarray}
Thus, equation (\ref{eq:mm_eq_mm}) holds, and $\mathfrak{P}$ is stationary.

Suppose that $(\mathfrak{m},\xi)\in\mathcal{S}_3$.
In the particular case $\alpha=1$, $\lambda=1$ and $c=0$, the
stationarity of $\mathfrak{P}$ was proved in Theorem~2 of~\cite
{kabluchkoetal}. The general case follows by a straightforward
application of affine transformations.
\end{pf}

\subsection{Two lemmas}
The next two lemmas are standard. We include their proofs only for completeness.
%
\begin{lemma}\label{lem:crit_stat_incr}
The process $W(t):=\xi(t)-\mu(t)$ has stationary increments iff for all
$t_1,t_2,h\in\mathbb{R}^d$,
%
\begin{equation}\label{eq:lem:gamma_stat}
\gamma(t_1,t_2)=\gamma(t_1+h,t_2+h).
\end{equation}
\end{lemma}
\begin{pf}
We prove only sufficiency since the necessity is evident. So, assume
that~(\ref{eq:lem:gamma_stat}) holds. Let $W_h(t)=W(t+h)-W(h)$. We have
\begin{eqnarray*}
&&\operatorname{Cov}(W_h(t_1), W_h(t_2))\\
&&\qquad=r(t_1+h,t_2+h)+r(h,h)-r(h,t_1+h)-r(h,t_2+h)\\
&&\qquad=- \bigl(\gamma(t_1+h,t_2+h)-\gamma(h,t_1+h)-\gamma(h, t_2+h) \bigr)/2\\
&&\qquad=- \bigl(\gamma(t_1,t_2)-\gamma(0,t_1)-\gamma(0,t_2) \bigr)/2,
\end{eqnarray*}
where the second equality follows from~(\ref{eq:r_sigma_gamma}) and
$\gamma(h,h)=0$, and the last equality is a consequence of~(\ref
{eq:lem:gamma_stat}).
Hence, the law of the process $\{W_h(t), t\in\mathbb{R}^d\}$ is
independent of
$h$, which proves the lemma.
\end{pf}

\begin{lemma}\label{lem:additive}
Let $g\dvtx \mathbb{R}^d\to\mathbb{R}$ be a function satisfying
%
\begin{equation}\label{eq:g_additive}
g(t_2+h)-g(t_1+h)=g(t_2)-g(t_1)
\end{equation}
for all $t_1,t_2,h\in\mathbb{R}^d$.
Then the following statements hold:
\begin{enumerate}
\item The function $f(t):=g(t)-g(0)$ is additive.
\item Either $g \equiv\mbox{ \textup{const}}$ or the set of values of $g$ is
dense in $\mathbb{R}$.
\end{enumerate}
\end{lemma}
\begin{pf}
Inserting $t_2:=s_1$, $h:=s_2$, $t_1:=0$ into~(\ref{eq:g_additive})
yields $f(s_1+s_2)=f(s_1)+f(s_2)$ and proves the first part of the
lemma. To prove the second part, assume that $g$ is not constant, which
means that there is $t$ with $f(t)\neq0$. A standard inductive
argument using the additivity of $f$ gives $f(qt)=qf(t)$ for every
rational number $q$. This implies that the set of values of the
function $f$, and hence also the set of values of $g$, is dense in
$\mathbb{R}$.
\end{pf}

\subsection{\texorpdfstring{Proof of Theorem~\protect\ref{theo:main_gen}:
Identifying the driving process $\xi$}{Proof of Theorem~1.1:
Identifying the driving process $\xi$}}

In Section~\ref{subsec:easy_part} we have shown that $\mathcal
{S}_1\cup\mathcal{S}
_2\cup
\mathcal{S}_3\subset\mathcal{S}$.
Here we prove the more difficult converse inclusion under an additional
assumption on the measure $\mathfrak{m}$. This is stated in the
following proposition.
\begin{proposition}\label{prop:proof_main_2}
Let $\mathfrak{m}$ be a measure of the form $\mathfrak{m}=\alpha
\mathfrak{e}_{\lambda}+\beta
\mathfrak{e}
_0$ for some $\alpha\geq0$, $\beta\geq0$, $\lambda\neq0$,
and let $\{\xi(t), t\in\mathbb{R}^d\}$ be a Gaussian process. Assume that
$\mathfrak{P}
=\mathit{GS}(\mathfrak{m},\xi)$ is stationary. Then $(\mathfrak{m},\xi)\in
\mathcal{S}_1\cup\mathcal{S}
_2\cup\mathcal{S}
_3$, where $\mathcal{S}_1,\mathcal{S}_2,\mathcal{S}_3$ are as in
Theorem~\ref{theo:main_gen}.
\end{proposition}

We will need some technical lemmas on measures which are obtained by
taking mixtures of diagonally shifted and exponentially weighted
bivariate normal laws.
%
\begin{lemma}\label{lem:proof_main}
Let $\mathfrak{n}$ be the law of a bivariate Gaussian vector
$(X_1,X_2)$ with
$\mathbb EX_i=\mu_i$, $\operatorname{Var}X_i=\sigma_i^2$ for $i=1,2$
and $\operatorname{Var}
(X_1-X_2)=\gamma$.
Let $\mathfrak{l}$ be a measure on $\mathbb{R}^2$ defined for some
$\kappa\in\mathbb{R}$ by
%
\begin{equation}\label{eq:def_lger}
\mathfrak{l}(B)=\int_{\mathbb{R}} e^{-\kappa z}\mathfrak{n}(B-z)\,dz,\qquad
B\in\mathcal{B}(\mathbb{R}^2).
\end{equation}
Then there is a measure $\mathfrak{l}^{(\kappa)}$ concentrated on the line
$\{
(x_1,x_2)\in\mathbb{R}^2\dvtx  x_1=0\}$ such that the following
representation holds:
%
\begin{equation}\label{eq:lger}
\mathfrak{l}(B)=\int_{\mathbb{R}} e^{-\kappa z} \mathfrak
{l}^{(\kappa)}(B-z)\,dz,\qquad B\in\mathcal{B}
(\mathbb{R}^2).
\end{equation}
The Laplace transform of $\mathfrak{l}^{(\kappa)}$, defined as $\psi
^{(\kappa
)}(u)=\int_{\mathbb{R}^2}e^{ux_2}\mathfrak{l}^{(\kappa
)}(dx_1,dx_2)$, is given by
%
\begin{equation}\label{eq:psi_star}
\hspace*{15pt}\psi^{(\kappa)}(u)=\exp\bigl\{(\kappa-u)\bigl(\mu_1+\tfrac12\kappa\sigma_1^2\bigr)
+u\bigl(\mu_2+\tfrac12 \kappa\sigma_2^2\bigr)+\tfrac12 u(u-\kappa)\gamma\bigr\}.
\end{equation}
\end{lemma}
%
\begin{remark}
Equation~(\ref{eq:psi_star}) shows that the measure $\mathfrak
{l}^{(\kappa)}$
is a multiple of a two-dimensional Gaussian measure. 
\end{remark}
%
\begin{remark}
If the Gaussian measure $\mathfrak{n}$ has a density, then it is
possible to
compute the density of $\mathfrak{l}$ directly from its definition,
equation~(\ref{eq:def_lger}). However, since $\mathfrak{n}$ (and also
$\mathfrak{l}$)
may fail to have a density, we use a somewhat more complicated
representation of $\mathfrak{l}$ as an exponentially weighted shift of the
essentially one-dimensional measure $\mathfrak{l}^{(\kappa)}$ given
in~(\ref{eq:lger}).
\end{remark}
\begin{pf*}{Proof of Lemma~\ref{lem:proof_main}}
Define
%
\begin{equation}\label{eq:l_star_def}
\mathfrak{l}^{(\kappa)}(B)=\int_{\mathbb{R}^2} e^{\kappa x_1}
1_{B}(0,x_2-x_1)\mathfrak{n}
(dx_1,dx_2),\qquad B\in\mathcal{B}(\mathbb{R}^2).
\end{equation}
By construction, the measure $\mathfrak{l}^{(\kappa)}$ is
concentrated on the
line $\{(x_1,x_2)\in\mathbb{R}^2\dvtx x_1=0\}$. Using transformations
similar to
those in~\cite{kabluchkoetal} (see the proof of Proposition~6
therein), we obtain
%
\begin{eqnarray}\label{eq:143}
\mathfrak{l}(B)
&=&\int_{\mathbb{R}}\int_{\mathbb{R}^2} e^{-\kappa z}
1_{B-z}(x_1,x_2)\mathfrak{n}(dx_1,dx_2)\,dz
\nonumber\\
&=&\int_{\mathbb{R}}\int_{\mathbb{R}^2} e^{-\kappa(z+x_1)}
e^{\kappa x_1}
1_{B-(z+x_1)}(0,x_2-x_1)\mathfrak{n}(dx_1,dx_2)\,dz \\
&=&\int_{\mathbb{R}}\int_{\mathbb{R}^2} e^{-\kappa w} e^{\kappa x_1}
1_{B-w}(0,x_2-x_1)\mathfrak{n}(dx_1,dx_2)\,dw. \nonumber
\end{eqnarray}
Applying~(\ref{eq:l_star_def}) to the right-hand side of the above
equation, we obtain~(\ref{eq:lger}).

Now we compute $\psi^{(\kappa)}(u)$, the Laplace transform of
$\mathfrak{l}
^{(\kappa)}$.
The Laplace transform of $\mathfrak{n}$ is defined as
\[
\psi(u_1,u_2)=\int_{\mathbb{R}^2}e^{u_1x_1+u_2x_2}\mathfrak{n}(dx_1,dx_2).
\]
By a two-dimensional analogue of~(\ref{eq:laplace_gauss}), $\psi
(u_1,u_2)$ is given by
%
\begin{equation}\label{eq:psi_laplace}
\psi(u_1,u_2)
=\exp\bigl\{\mu_1u_1+\mu_2u_2+\tfrac12(\sigma_1^2 u_1^2+2ru_1u_2+\sigma_2^2u_2^2) \bigr\},
\end{equation}
where $r=\operatorname{Cov}(X_1,X_2)=(\sigma_1^2+\sigma_2^2-\gamma)/2$.
It follows from~(\ref{eq:l_star_def}) that
\[
\psi^{(\kappa)}(u)=\int_{\mathbb{R}^2} e^{\kappa x_1}
e^{u(x_2-x_1)} \mathfrak{n}
(dx_1,dx_2)=\psi(\kappa-u,u).
\]
The above equation and~(\ref{eq:psi_laplace}) yield~(\ref{eq:psi_star})
after an elementary calculation.
\end{pf*}

%
\begin{lemma}\label{lem:proof_main_uniq}
Fix $\kappa\neq0$. Let $\mathfrak{l}$ be a Radon measure on $\mathbb
{R}^2$ admitting
a decomposition
%
\begin{equation}\label{eq:lger_uniq}
\mathfrak{l}(B)=\int_{\mathbb{R}} e^{-\kappa z} \mathfrak
{l}^{(\kappa)}(B-z)\,dz+\int_{\mathbb{R}
} \mathfrak{l}
^{(0)}(B-z)\,dz,\qquad
B\in\mathcal{B}(\mathbb{R}^2),
\end{equation}
where $\mathfrak{l}^{(\kappa)}$ and $\mathfrak{l}^{(0)}$ are
measures concentrated on
the line $\{(x_1,x_2)\in\mathbb{R}^2:x_1=0\}$. Then the measures
$\mathfrak{l}
^{(\kappa
)}$ and $\mathfrak{l}^{(0)}$ are determined uniquely.
\end{lemma}
\begin{pf}
Fix some bounded Borel set $A\subset\{0\}\times\mathbb{R}$. For
$x>0$, let
$B_x$ be a subset of $\mathbb{R}^2$ defined by $B_x=\bigcup_{y\in[0,x]}(A+y)$.
Then~(\ref{eq:lger_uniq}) implies that
\[
\mathfrak{l}(B_x)= \biggl(\int_{0}^x e^{-\kappa z}\,dz \biggr) \mathfrak
{l}^{(\kappa)}(A) +x
\mathfrak{l}^{(0)}(A).
\]
The above is valid for every $x>0$, and so, $\mathfrak{l}^{(\kappa
)}(A)$ and
$\mathfrak{l}^{(0)}(A)$ are determined uniquely.
\end{pf}

\begin{pf*}{Proof of Proposition~\ref{prop:proof_main_2}}
We start by proving three claims about the expectation $\mu(\cdot)$,
the variance $\sigma^2(\cdot)$ and the incremental variance $\gamma
(\cdot,\cdot)$ under various assumptions on $\alpha,\beta,\lambda$.
%
\begin{claim}\label{claim:1}
Assume that $\alpha>0$. Then for all $t_1,t_2\in\mathbb{R}^d$,
%
\begin{equation}
\mu(t_2)-\mu(t_1)=-\frac{\lambda}{2} \bigl(\sigma^2(t_2)-\sigma^2(t_1) \bigr).
\label{eq:mu_sigma_proof}
\end{equation}
\end{claim}
\begin{pf}
The measure $\mathfrak{m}_t=\mathfrak{m}*\mathfrak{n}_t$ has a
density given by the convolution formula
\[
\frac{\mathfrak{m}_t(dx)}{dx}
=\int_{\mathbb{R}} \bigl(\alpha e^{-\lambda(x-y)}+\beta\bigr)\mathfrak{n}_t(dy)
=\alpha e^{-\lambda x} \int_{\mathbb{R}} e^{\lambda y} \mathfrak
{n}_t(dy)+ \beta.
\]
Applying~(\ref{eq:laplace_gauss}) to the first term on right-hand side,
we obtain
%
\begin{equation}\label{eq:dens344}
\frac{\mathfrak{m}_t(dx)}{dx}=\alpha e^{-\lambda x} \exp\biggl\{\mu(t)
\lambda
+\frac
12 \sigma^2(t)\lambda^2 \biggr\}+\beta.
\end{equation}
%
By stationarity of $\mathfrak{P}$, we must have $\mathfrak
{m}_{t_1}=\mathfrak{m}_{t_2}$ for every
$t_1,t_2\in\mathbb{R}^d$. This leads to~(\ref{eq:mu_sigma_proof}).
\end{pf}


Let us turn to the ``two-dimensional distributions'' of $\mathfrak
{P}$. Take
$t_1,t_2\in\mathbb{R}^d$ and recall that $\mathfrak{P}_{t_1,t_2}=\{(V_i(t_1),
V_i(t_2)), i\in\mathbb{N}\}$ is a Poisson point process on $\mathbb
{R}^2$. By~(\ref{eq:mm_t_R}), its intensity measure $\mathfrak
{m}_{t_1,t_2}$ is given for $B\in
\mathcal{B}(\mathbb{R}^2)$ by
\begin{eqnarray}\label{eq:mt1t2}
\mathfrak{m}_{t_1,t_2}(B)
&=&\int_{\mathbb{R}} (\alpha e^{-\lambda x}+\beta)\mathfrak
{n}_{t_1,t_2}(B-x)\,dx\nonumber
\\[-8pt]\\[-8pt]
&=&\alpha\int_{\mathbb{R}} e^{-\lambda x}\mathfrak
{n}_{t_1,t_2}(B-x)\,dx+ \beta\int
_{\mathbb{R}
}\mathfrak{n}_{t_1,t_2}(B-x)\,dx. \nonumber
\end{eqnarray}
Applying Lemma~\ref{lem:proof_main} twice with $\kappa=\lambda$,
$\mathfrak{n}
=\mathfrak{n}_{t_1,t_2}$ and $\kappa=0$, $\mathfrak{n}=\mathfrak
{n}_{t_1,t_2}$, we obtain two
measures on $\mathbb{R}^2$, called $\mathfrak{m}_{t_1,t_2}^{(\lambda
)}$ and $\mathfrak{m}
^{(0)}_{t_1,t_2}$, which are concentrated on the line $\{(x_1,x_2)\in
\mathbb{R}
^2\dvtx  x_1=0\}$ and have the property that for each Borel set $B\subset
\mathbb{R}^2$,
%
\begin{equation}\label{eq:mt1t2_2}
\mathfrak{m}_{t_1,t_2}(B)=\alpha\int_{\mathbb{R}} e^{-\lambda x}
\mathfrak{m}
_{t_1,t_2}^{(\lambda
)}(B-x)\,dx + \beta\int_{\mathbb{R}} \mathfrak{m}_{t_1,t_2}^{(0)}(B-x)\,dx.
\end{equation}
%
%
\begin{claim}\label{claim:2}
Assume that $\alpha>0$. Then for all $t_1,t_2,h\in\mathbb{R}^d$,
%
\begin{equation}
\gamma(t_1,t_2)=\gamma(t_1+h,t_2+h). \label{eq:gamma_stat}
\end{equation}
\end{claim}
\begin{pf}
By stationarity, $\mathfrak{m}_{t_1,t_2}=\mathfrak{m}_{t_1+h,t_2+h}$
for all
$t_1,t_2,h\in
\mathbb{R}^d$. Applying Lemma~\ref{lem:proof_main_uniq} to the
decomposition~(\ref{eq:mt1t2_2}), we obtain
\[
\mathfrak{m}_{t_1,t_2}^{(\lambda)}=\mathfrak
{m}_{t_1+h,t_2+h}^{(\lambda)}.
\]
Recall that\vspace*{1pt} the measures $\mathfrak{m}_{t_1,t_2}^{(\lambda)}$ and
$\mathfrak{m}
_{t_1+h,t_2+h}^{(\lambda)}$ were constructed by means of Lemma~\ref
{lem:proof_main} and thus have Laplace transforms given by the
right-hand side of~(\ref{eq:psi_star}). So, we obtain that the
expression (considered as a polynomial in $u$)
\[
(\lambda-u) \biggl(\mu(t_1)+\frac{\lambda}2 \sigma^2(t_1) \biggr)
+u \biggl(\mu(t_2)+\frac{\lambda}2 \sigma^2(t_2) \biggr)+\frac12u(u-\lambda)\gamma(t_1,t_2)
\]
does not change if we replace $t_1,t_2$ by $t_1+h,t_2+h$. Taking into
account that by Claim~\ref{claim:1},
\[
\mu(t_i)+\frac{\lambda}2\sigma^2(t_i)=\mu(t_i+h)+\frac{\lambda}2
\sigma^2(t_i+h),\qquad i=1,2,
\]
we arrive at~(\ref{eq:gamma_stat}).
\end{pf}
%
%
\begin{claim}\label{claim:3}
Assume that $\beta>0$. Then for all $t_1,t_2,h\in\mathbb{R}^d$,
%
\begin{equation}\label{eq:mu_stat}
\mu(t_2)-\mu(t_1)=\mu(t_2+h)-\mu(t_1+h)
\end{equation}
and
%
\begin{equation}\label{eq:gamma_stat_claim:3}
\gamma(t_1,t_2)=\gamma(t_1+h,t_2+h).
\end{equation}
\end{claim}
\begin{pf}
It follows from the decomposition~(\ref{eq:mt1t2_2}) and Lemma~\ref
{lem:proof_main_uniq} that
\[
\mathfrak{m}_{t_1,t_2}^{(0)}=\mathfrak{m}_{t_1+h,t_2+h}^{(0)}.
\]
Using the formula\vspace*{1pt} for the Laplace transform of $\mathfrak{m}_{t_1,t_2}^{(0)}$
and $\mathfrak{m}_{t_1+h,t_2+h}^{(0)}$ given in~(\ref{eq:psi_star}),
we obtain
that the expression (considered as a quadratic polynomial in $u$)
\[
u\bigl(\mu(t_2)-\mu(t_1)\bigr)+\tfrac12 \gamma(t_1,t_2)u^2 
\]
remains unchanged if we replace $t_1,t_2$ by $t_1+h,t_2+h$. This
yields~(\ref{eq:mu_stat}) and~(\ref{eq:gamma_stat_claim:3}).
\end{pf}

Now we are ready to complete the proof of Proposition~\ref
{prop:proof_main_2}. We distinguish three cases.
\begin{case}
Assume that $\alpha>0$ and $\beta>0$. We show that in this case,
$(\mathfrak{m}
,\xi)\in\mathcal{S}_1$. Combining Claims~\ref{claim:1} and \ref{claim:3},
we obtain
\[
\sigma^2(t_2)-\sigma^2(t_1)=\sigma^2(t_2+h)-\sigma^2(t_1+h).
\]
Since $\sigma^2(t)\geq0$, it follows from part~2 of Lemma~\ref
{lem:additive} that $\sigma^2(t)$ is a constant function. By
Claim~\ref
{claim:1}, $\mu(t)$ is constant as well. Finally, by~(\ref
{eq:r_sigma_gamma}) and Claim~\ref{claim:2},
\begin{eqnarray*}
r(t_1+h,t_2+h)
&=&\tfrac12 \bigl(\sigma^2(t_1+h)+\sigma^2(t_2+h)-\gamma(t_1+h,t_2+h) \bigr)\\
&=&\tfrac12 \bigl(2\sigma^2(0)-\gamma(t_1,t_2) \bigr)\\
&=&r(t_1,t_2).
\end{eqnarray*}
This implies that the Gaussian process $W(t):=\xi(t)-\mu(t)$ is
stationary. Hence, $(\mathfrak{m},\xi)\in\mathcal{S}_1$.
\end{case}

\begin{case}
Assume that $\alpha=0$ and $\beta>0$. We show that in this case,
$(\mathfrak{m}
,\xi)\in\mathcal{S}_2$. First of all, note that in this case,
$\mathfrak{m}$ is a
multiple of $\mathfrak{e}_0$. By equation~(\ref
{eq:gamma_stat_claim:3}) of
Claim~\ref{claim:3} and Lemma~\ref{lem:crit_stat_incr}, the process
$W(t):=\xi(t)-\mu(t)$ has stationary increments. Further, the function
$f(t):=\mu(t)-\mu(0)$ is additive by equation~(\ref{eq:mu_stat}) of
Claim~\ref{claim:3} and part~1 of Lemma~\ref{lem:additive}. So, we
obtain a decomposition $\xi(t)=W(t)+f(t)+\mu(0)$ implying that
$(\mathfrak{m}
,\xi
)\in\mathcal{S}_2$.
\end{case}

\begin{case}
Assume that $\alpha>0$ and $\beta=0$. We show that in this case,
$(\mathfrak{m}
,\xi)\in\mathcal{S}_3$. First, we have $\mathfrak{m}=\alpha
\mathfrak{e}_{\lambda}$. Second,
Claim~\ref{claim:2} and Lemma~\ref{lem:crit_stat_incr} show that the
process $W(t):=\xi(t)-\mu(t)$ has stationary increments. It follows
from Claim~\ref{claim:1} that
\[
\mu(t)=-\lambda\sigma^2(t)/2+\mu(0)+\lambda\sigma^2(0)/2=-\lambda
\sigma^2(t)/2+c,
\]
where $c=\mu(0)+\lambda\sigma^2(0)/2$. Hence, $(\mathfrak{m},\xi
)\in\mathcal{S}_3$.
\end{case}

The proof of Proposition~\ref{prop:proof_main_2} is complete.
\end{pf*}

\subsection{Lemmas on convolution equations}\label{subsec:lem_conv}
In this section we collect several auxiliary lemmas on solutions of
convolution equations. These equations will arise in Section~\ref
{subsec:ident_m} when dealing with one-dimensional distributions of
Gaussian systems.
The proofs are based on explicit calculations with Laplace transforms
and on the result of Deny~\cite{deny59}.
%
\begin{lemma}\label{lem:conv_eq_3}
Let $\mathfrak{n}_0=\mathfrak{N}(\mu_0,\sigma_0^2)$ be a Gaussian
measure on $\mathbb{R}$. Let
$\mathfrak{m}_1,\mathfrak{m}_2$ be two measures satisfying~(\ref
{eq:mm_cond_fin}) such that
%
\begin{equation}\label{eq:conv_eq_2}
\mathfrak{m}_1*\mathfrak{n}_0=\mathfrak{m}_2*\mathfrak{n}_0.
\end{equation}
Then $\mathfrak{m}_1=\mathfrak{m}_2$.
\end{lemma}
\begin{pf}
We assume that $\sigma_0^2>0$, since otherwise, the statement of the
lemma is trivial. The density of the measure $\mathfrak{m}_i*\mathfrak
{n}_0$, $i=1,2$, is
given by the convolution formula
%
\begin{eqnarray}\label{eq:mm_conv_nn_i_dens}
\frac{(\mathfrak{m}_i*\mathfrak{n}_0)(dx)}{dx}
&=&\frac{1}{\sqrt{2\pi}\sigma_0} \int_{\mathbb{R}}e^{-\afrac{(x-y-\mu_0)^2}{2\sigma_0^2}}\mathfrak{m}_i(dy)\nonumber\\
&=&\frac{1}{\sqrt{2\pi}\sigma_0} e^{-\afrac{x^2}{2\sigma_0^2}}
e^{\sfrac{x\mu_0}{\sigma_0^2}}
\\
&&{}\times\int_{\mathbb{R}}e^{\sfrac{xy}{\sigma_0^2}}e^{-\afrac{(y+\mu_0)^2}{2\sigma_0^2}} \mathfrak{m}_i(dy). \nonumber
\end{eqnarray}
Define new measures $\mathfrak{m}_1'$ and $\mathfrak{m}_2'$ by
%
\begin{equation}\label{eq:density_mm_i_prime_def}
\frac{\mathfrak{m}_i'(dy)}{\mathfrak{m}_i(dy)}=e^{-\afrac{(y+\mu_0)^2}{2\sigma_0^2}},\qquad i=1,2.
\end{equation}
Let $\varphi_{\mathfrak{m}_i'}(x)=\int_{\mathbb{R}}e^{xy}\mathfrak
{m}_i'(dy)$, $i=1,2$, be the
Laplace transforms of $\mathfrak{m}_1'$ and $\mathfrak{m}_2'$. Note
that by~(\ref{eq:mm_cond_fin}), $\varphi_{\mathfrak{m}_1'}(x)$ and
$\varphi_{\mathfrak{m}_2'}(x)$ are
finite for all $x\in\mathbb{R}$. We may rewrite~(\ref{eq:mm_conv_nn_i_dens})
as follows:
\begin{eqnarray}\label{eq:mm_conv_nn_i_dens_cont}
\frac{(\mathfrak{m}_i*\mathfrak{n}_0)(dx)}{dx}
&=&
\frac{1}{\sqrt{2\pi}\sigma_0} e^{-\afrac{x^2}{2\sigma_0^2}}e^{\sfrac{x\mu_0}{\sigma_0^2}}
\int_{\mathbb{R}}e^{\sfrac{xy}{\sigma_0^2}}\mathfrak{m}_i'(dy)\nonumber  \\[-8pt]\\[-8pt]
&=&
\frac{1}{\sqrt{2\pi}\sigma_0} e^{-\afrac{x^2}{2\sigma_0^2}}
e^{\sfrac{x\mu_0}{\sigma_0^2}} \varphi_{\mathfrak{m}_i'} \biggl(\frac{x}{\sigma_0^2} \biggr).
\nonumber
\end{eqnarray}
By~(\ref{eq:conv_eq_2}), the densities of the measures $\mathfrak
{m}_1*\mathfrak{n}_0$
and $\mathfrak{m}_2*\mathfrak{n}_0$ must be equal.
Taking into account~(\ref{eq:mm_conv_nn_i_dens_cont}), this yields
\[
\varphi_{\mathfrak{m}_1'}(x)= \varphi_{\mathfrak{m}_2'}(x) \qquad\forall x\in\mathbb{R}.
\]
By the uniqueness of the Laplace transform, $\mathfrak{m}'_1=\mathfrak{m}'_2$.
Recalling~(\ref{eq:density_mm_i_prime_def}) yields that $\mathfrak
{m}_1=\mathfrak{m}_2$.
This proves the lemma.
\end{pf}
%
\begin{lemma}\label{lem:conv_eq_4}
Let $\mathfrak{n}_1=\mathfrak{N}(\mu_1,\sigma_1^2)$ and $\mathfrak
{n}_2=\mathfrak{N}(\mu_2,\sigma_2^2)$
be two Gaussian measures on $\mathbb{R}$ such that $\sigma_1^2\leq
\sigma_2^2$.
Let $\mathfrak{m}_1$ and $\mathfrak{m}_2$ be two measures
satisfying~(\ref{eq:mm_cond_fin}) such that
%
\begin{equation}\label{eq:conv_eq_4}
\mathfrak{m}_1*\mathfrak{n}_1=\mathfrak{m}_2*\mathfrak{n}_2.
\end{equation}
Then $\mathfrak{m}_1=\mathfrak{m}_2*\mathfrak{N}(\mu_2-\mu
_1,\sigma_2^2-\sigma_1^2)$.
\end{lemma}
\begin{pf}
We may rewrite~(\ref{eq:conv_eq_4}) as
\[
\mathfrak{m}_1*\mathfrak{n}_1= \bigl(\mathfrak{m}_2*\mathfrak{N}(\mu
_2-\mu_1,\sigma_2^2-\sigma_1^2)
\bigr)*\mathfrak{n}_1.
\]
The proof is completed by applying Lemma~\ref{lem:conv_eq_3}.
\end{pf}

\begin{lemma}\label{lem:conv_eq_1}
Let $\mathfrak{m}$ be a measure satisfying~(\ref{eq:mm_cond_fin}),
and let
$\mathfrak{n}
_0=\mathfrak{N}(\mu_0,\sigma_0^2)$ be a Gaussian measure such that
for some
$\alpha\geq0$, $\beta\geq0$, $\lambda\neq0$,
%
\begin{equation}\label{eq:conv_lem_1}
\mathfrak{m}*\mathfrak{n}_0=\alpha\mathfrak{e}_{\lambda}+\beta
\mathfrak{e}_0.
\end{equation}
Then $\mathfrak{m}=\alpha e^{-\lambda^2\sigma^2_0/2}e^{-\lambda\mu
_0}\mathfrak{e}
_{\lambda}+\beta\mathfrak{e}_0$.
\end{lemma}
\begin{pf}
Define a measure $\mathfrak{m}_1=\alpha e^{-\lambda^2\sigma
^2_0/2}e^{-\lambda
\mu
_0}\mathfrak{e}_{\lambda}+\beta\mathfrak{e}_0$. Then the density of
the measure $\mathfrak{m}
_1*\mathfrak{n}_0$ can be computed by means of the convolution formula:
\[
\frac{(\mathfrak{m}_1*\mathfrak{n}_0)(dx)}{dx}=\int_{\mathbb
{R}}\bigl(\alpha e^{-\lambda^2\sigma
^2_0/2}e^{-\lambda\mu_0}e^{-\lambda(x-y)}+\beta\bigr)\mathfrak
{n}_0(dy)=\alpha
e^{-\lambda x}+\beta,
\]
where the second equality follows from~(\ref{eq:laplace_gauss}). Hence,
\[
\mathfrak{m}*\mathfrak{n}_0=\mathfrak{m}_1*\mathfrak{n}_0.
\]
By Lemma~\ref{lem:conv_eq_3}, we have $\mathfrak{m}=\mathfrak{m}_1$.
The proof is complete.
\end{pf}
%
\begin{lemma}\label{lem:conv_eq_2}
Let $\mathfrak{n}_1=\mathfrak{N}(\mu_1,\sigma_1^2)$ and $\mathfrak
{n}_2=\mathfrak{N}(\mu_2,\sigma_2^2)$
be two Gaussian measures on $\mathbb{R}$ such that $\sigma_1^2\neq
\sigma_2^2$.
Let $\mathfrak{m}$ be a measure satisfying~(\ref{eq:mm_cond_fin})
such that
%
\begin{equation}\label{eq:conv_eq}
\mathfrak{m}*\mathfrak{n}_1=\mathfrak{m}*\mathfrak{n}_2.
\end{equation}
Then $\mathfrak{m}=\alpha\mathfrak{e}_{\lambda}+\beta\mathfrak
{e}_{0}$ for some $\alpha\geq0$,
$\beta\geq0$ and $\lambda\neq0$.
\end{lemma}
\begin{pf}
By symmetry, we may assume that $\sigma_1^2<\sigma_2^2$.
Then Lemma~\ref{lem:conv_eq_3} implies that
%
\begin{equation}\label{eq:conv_deny_prime}
\mathfrak{m}=\mathfrak{m}*\mathfrak{n}_0,
\end{equation}
where
$
\mathfrak{n}_0=\mathfrak{N}(\mu_2-\mu_1,\sigma_2^2-\sigma_1^2)$.
By Theorem~3$'$ of~\cite{deny59}, every solution $\mathfrak{m}$
of~(\ref{eq:conv_deny_prime}) can be represented as a mixture of exponentials;
that is, we may write
\[
\frac{\mathfrak{m}(dy)}{dy}=\int_{E}e^{-\lambda y} \rho(d\lambda),
\]
where $\rho$ is a finite Borel measure on the set $E=\{\lambda\in
\mathbb{R}\dvtx
\int_{\mathbb{R}}e^{\lambda x} \mathfrak{n}_0(dx)=1\}$.
Now, in our case the measure $\mathfrak{n}_0$ is Gaussian, and
so~(\ref{eq:laplace_gauss}) shows that $E$ consists of at most two
points. One
of them is always $0$, and the second is denoted by $\lambda$ (if $E=\{
0\}$, let $\lambda\neq0$ be arbitrary). Taking $\alpha=\rho(\{
\lambda\}
)$ and $\beta=\rho(\{0\})$, we obtain
$\mathfrak{m}=\alpha\mathfrak{e}_{\lambda}+\beta\mathfrak{e}_0$.
This completes the proof.
\end{pf}

\subsection{\texorpdfstring{Proof of Theorem~\protect\ref{theo:main_gen}:
Identifying the measure $\mathfrak{m}$}{Proof of Theorem~1.1:
Identifying the measure $\mathfrak{m}$}}
\label{subsec:ident_m}

In this section we complete the proof of the inclusion $\mathcal
{S}\subset\mathcal{S}
_1\cup\mathcal{S}_2\cup\mathcal{S}_3$. Let $(\mathfrak{m},\xi)$
be a pair generating a
stationary Gaussian system $\mathfrak{P}=\mathit{GS}(\mathfrak{m},\xi)$. Our
goal is to show that
%
\begin{equation}\label{eq:goal}
(\mathfrak{m},\xi)\in\mathcal{S}_1\cup\mathcal{S}_2\cup\mathcal{S}_3.
\end{equation}
The idea of the proof is to show, whenever possible, that the measure
$\mathfrak{m}$ is of the form $\alpha\mathfrak{e}_{\lambda}+\beta
\mathfrak{e}_{0}$ and then to
apply Proposition~\ref{prop:proof_main_2}. In all other cases, we will
prove that $(\mathfrak{m},\xi)\in\mathcal{S}_1$.

Assume for a moment that $\xi(0)=0$ and $\operatorname{Var}\xi
(t_0)>0$ for some
$t_0\in\mathbb{R}^d$. Under this restriction, the proof takes the following
particularly simple form.
By stationarity, we have $\mathfrak{m}_0=\mathfrak{m}_{t_0}$. Using
$\xi(0)=0$, this can
be written as $\mathfrak{m}=\mathfrak{m}*\mathfrak{n}_{t_0}$.
Applying to this convolution
equation the result of Deny~\cite{deny59} as in the proof of
Lemma~\ref
{lem:conv_eq_2}, we conclude that $\mathfrak{m}$ must be of the form
$\alpha
\mathfrak{e}
_{\lambda}+\beta\mathfrak{e}_0$. Hence, Proposition~\ref
{prop:proof_main_2} is
applicable and~(\ref{eq:goal}) is proved.

Let us now consider Theorem~\ref{theo:main_gen} in its full generality.
We will distinguish between different cases.
%
\begin{casee}
Assume that the function $\sigma^2$ is not constant. So, there are
$t_1,t_2\in\mathbb{R}^d$ such that
%
\begin{equation}\label{eq:sigma_not_const}
\sigma^2(t_1)\neq\sigma^2(t_2).
\end{equation}
By stationarity of $\mathfrak{P}$, we must have $\mathfrak
{m}_{t_1}=\mathfrak{m}_{t_2}$ and hence,
\[
\mathfrak{m}*\mathfrak{N}(\mu(t_1),\sigma^2(t_1))=\mathfrak
{m}*\mathfrak{N}(\mu(t_2),\sigma^2(t_2)).
\]
Then Lemma~\ref{lem:conv_eq_2}, which is applicable in view of~(\ref
{eq:sigma_not_const}), implies that $\mathfrak{m}=\alpha\mathfrak
{e}_{\lambda}+\beta
\mathfrak{e}
_{0}$ for some $\alpha\geq0$, $\beta\geq0$, $\lambda\neq0$.
An application of Proposition~\ref{prop:proof_main_2} shows
that~(\ref{eq:goal}) holds.
\end{casee}

\begin{casee}
Assume that\vspace*{1pt} $\sigma^2(t)=\sigma^2\geq0$ is constant.
Take some $t_1,t_2\in\mathbb{R}^d$ and fix some $\vartheta\in[0,1]$.
Consider $\tilde\mathfrak{P}_{t_1,t_2}$, a point process on $\mathbb
{R}$ defined by
%
\begin{equation}\label{eq:def_P_wsp110}
\tilde\mathfrak{P}_{t_1,t_2}= \{U_i+ \vartheta\xi
_i(t_1)+(1-\vartheta)\xi
_i(t_2), i\in\mathbb{N}\},
\end{equation}
where the $U_i$'s and the $\xi_i$'s are as in Section~\ref{subsec:problem}.
Recalling from~(\ref{eq:def_v}) that $V_i(t)=U_i+\xi_i(t)$, we may
rewrite~(\ref{eq:def_P_wsp110}) as
%
\begin{equation}\label{eq:def_P_wsp111}
\tilde\mathfrak{P}_{t_1,t_2}= \{\vartheta V_i(t_1)+(1-\vartheta) V_i(t_2),
i\in
\mathbb{N}\}.
\end{equation}
By Proposition~3.8 of~\cite{resnickbook87}, $\tilde\mathfrak
{P}_{t_1,t_2}$ is
a Poisson point process whose intensity measure $\tilde\mathfrak{m}_{t_1,t_2}$
is given by the formula
\[
\tilde\mathfrak{m}_{t_1,t_2}=\mathfrak{m}*\mathfrak{N}(\tilde\mu
(t_1,t_2),\tilde\sigma
^2(t_1,t_2)),
\]
where
%
\begin{eqnarray}
\tilde\mu(t_1,t_2)
&=\vartheta\mu(t_1)+(1-\vartheta)\mu(t_2)\label{eq:mu_tilde}
\end{eqnarray}
and
%
\begin{eqnarray}
\tilde\sigma^2(t_1,t_2)
=\bigl(\vartheta^2+(1-\vartheta)^2\bigr)\sigma^2+2\vartheta(1-\vartheta
)r(t_1, t_2).
\label{eq:sigma_tilde}
\end{eqnarray}
The stationarity of the particle system $\mathfrak{P}$ together with
representation~(\ref{eq:def_P_wsp111}) implies that for every
$t_1,t_2,h\in\mathbb{R}^d$, the point processes $\tilde\mathfrak
{P}_{t_1,t_2}$ and
$\tilde\mathfrak{P}_{t_1+h, t_2+h}$ must have the same law. Hence,
$\tilde
\mathfrak{m}
_{t_1,t_2}=\tilde\mathfrak{m}_{t_1+h, t_2+h}$ and consequently,
%
\begin{equation}\label{eq:mm_conv_sigma_tilde}
\hspace*{28pt}\mathfrak{m}*\mathfrak{N}(\tilde\mu(t_1,t_2),\tilde\sigma
^2(t_1,t_2))=\mathfrak{m}*\mathfrak{N}
\bigl(\tilde
\mu(t_1+h,t_2+h),\tilde\sigma^2(t_1+h,t_2+h)\bigr).
\end{equation}
The proof will be completed after we have considered two subcases.
\begin{subcase}
Assume that for some $t_1,t_2,h\in\mathbb{R}^d$,
%
\begin{equation}\label{eq:r_not_stationary}
r(t_1, t_2)\neq r(t_1+h, t_2+h).
\end{equation}
Take $\vartheta=1/2$ in the definition of the point process $\tilde
\mathfrak{P}
_{t_1,t_2}$. Then~(\ref{eq:sigma_tilde}) and~(\ref
{eq:r_not_stationary}) imply that
\[
\tilde\sigma^2(t_1,t_2)\neq\tilde\sigma^2(t_1+h,t_2+h).
\]
By Lemma~\ref{lem:conv_eq_2}, applied to~(\ref
{eq:mm_conv_sigma_tilde}), the measure $\mathfrak{m}$ is of the form
$\alpha\mathfrak{e}
_{\lambda}+\beta\mathfrak{e}_{0}$ for some $\alpha\geq0$, $\beta
\geq0$,
$\lambda\neq0$. An application of Proposition~\ref{prop:proof_main_2}
shows that~(\ref{eq:goal}) holds.
\end{subcase}
%
\begin{subcase}
Assume that for all $t_1,t_2,h\in\mathbb{R}^d$,
%
\begin{equation}\label{eq:r_stationary}
r(t_1, t_2)=r(t_1+h, t_2+h).
\end{equation}
%
This implies that the process $W(t):=\xi(t)-\mu(t)$ is stationary.

If the function $\mu$ is constant, then $(\mathfrak{m},\xi)\in
\mathcal{S}_1$.
Therefore, let us assume that $\mu$ is not constant. We will show that
this implies that $\mathfrak{m}$ is a multiple of the Lebesgue measure.
Let
\[
G=\{g\in\mathbb{R}\dvtx  \mathfrak{m}*\delta_{g}=\mathfrak{m}\}
\]
be the set of ``periods'' of $\mathfrak{m}$, where $\delta_g$ is the Dirac
measure concentrated at $g$. Clearly, $G$ is an additive subgroup of
$\mathbb{R}$.

By stationarity of $\mathfrak{P}$, we have $\mathfrak
{m}_{t_1}=\mathfrak{m}_{t_2}$ for every
$t_1,t_2\in\mathbb{R}^d$. Equivalently,
\[
\mathfrak{m}*\mathfrak{N}(\mu(t_1),\sigma^2)=\mathfrak
{m}*\mathfrak{N}(\mu(t_2),\sigma^2).
\]
By Lemma~\ref{lem:conv_eq_4}, this implies that
%
\begin{equation}\label{eq:mu_t1_mu_t2_G}
\mu(t_1)-\mu(t_2)\in G \qquad\forall t_1,t_2\in\mathbb{R}^d.
\end{equation}

Since $\mu$ is assumed to be nonconstant, equation~(\ref
{eq:mu_t1_mu_t2_G}) implies that $G\neq\{0\}$, which means that
$\mathfrak{m}$
has a nontrivial period. Of course, this is not sufficient to conclude
that $\mathfrak{m}$ is a multiple of the Lebesgue measure, and so, let
us use
the stationarity of the two-dimensional distributions of $\mathfrak{P}$.
Recalling~(\ref{eq:sigma_tilde}) and taking into account~(\ref
{eq:r_stationary}), we obtain that for every $t_1,t_2,h\in\mathbb{R}^d$,
\[
\tilde\sigma^2(t_1,t_2)=\tilde\sigma^2(t_1+h,t_2+h).
\]
Applying Lemma~\ref{lem:conv_eq_4} to~(\ref{eq:mm_conv_sigma_tilde}),
we obtain
\[
\tilde\mu(t_1,t_2)-\tilde\mu(t_1+h,t_2+h)\in G \qquad\forall t_1,t_2,h\in
\mathbb{R}^d.
\]
Recalling a formula for $\tilde\mu$ given in~(\ref{eq:mu_tilde}), we
arrive at
\[
\vartheta\cdot\bigl(\mu(t_1)-\mu(t_2)-\mu(t_1+h)+\mu(t_2+h)\bigr)+\bigl(\mu
(t_2)-\mu
(t_2+h)\bigr) \in G.
\]
Note that this is valid for every $\vartheta\in[0,1]$. Assume that in
the above expression, $\vartheta$ appears with a nonzero coefficient
for some $t_1,t_2,h$. Then $G$ contains a nontrivial interval, and so,
we must have $G=\mathbb{R}$. In other words, the measure $\mathfrak
{m}$ is translation
invariant. Since by~(\ref{eq:mm_cond_fin}), $\mathfrak{m}$ is finite
on bounded
intervals, this implies that $\mathfrak{m}$ is a multiple of the
Lebesgue measure.

So, let us assume that for every $t_1,t_2,h\in\mathbb{R}^d$,
%
\begin{equation}\label{eq:mu_additive}
\mu(t_1)-\mu(t_2)=\mu(t_1+h)-\mu(t_2+h).
\end{equation}
Recall also that we assume that $\mu$ is nonconstant. Hence, by part~2
of Lemma~\ref{lem:additive}, the set of values of the function $\mu$ is
dense in $\mathbb{R}$. By~(\ref{eq:mu_t1_mu_t2_G}), the group $G$
must be dense
in $\mathbb{R}$.

We claim that in fact, $G=\mathbb{R}$. To prove this, we need to show
that $G$
is closed. First of all, the measure $\mathfrak{m}$ is atomless, since
if it
would have an atom, then the invariance under $G$ would imply that
$\mathfrak{m}
$ has a dense set of atoms of equal mass, which would contradict~(\ref
{eq:mm_cond_fin}). Now, let $g_1,g_2,\ldots$ be a sequence in $G$
converging to some $g\in\mathbb{R}$. We claim that $g\in G$. Indeed,
for every
interval $[a,b]\subset\mathbb{R}$, we have
\[
\mathfrak{m}([a-g,b-g])=\lim_{i\to\infty}\mathfrak
{m}([a-g_i,b-g_i])=\lim_{i\to
\infty}\mathfrak{m}
([a,b])=\mathfrak{m}([a,b]),
\]
where the first equality holds since $\mathfrak{m}$ is atomless, and
the second
equality follows from $g_i\in G$. This proves that $g\in G$. Therefore,
the group $G$, being dense and closed, must be equal to $\mathbb{R}$.

The fact that $G=\mathbb{R}$ means that the measure $\mathfrak{m}$ is
translation
invariant and thus, must be a multiple of the Lebesgue measure.
Therefore, we can apply Proposition~\ref{prop:proof_main_2} which shows
that~(\ref{eq:goal}) holds.
\end{subcase}
\end{casee}

The proof of Theorem~\ref{theo:main_gen} is complete.

\section{Pairs generating equal in law Gaussian systems}\label
{sec:class_complete}
In this section we give an answer to the following question: Given two
pairs $(\mathfrak{m}',\xi')$ and $(\mathfrak{m}'',\xi'')$ in
$\mathcal{S}$, determine whether
$\mathit{GS}(\mathfrak{m}',\xi')$ has the same law as $\mathit{GS}(\mathfrak{m}'',\xi
'')$. The next
proposition is a first step in this direction.
%
\begin{proposition}\label{prop:after_theo_main_gen}
The decomposition
$
\mathcal{S}=\mathcal{S}_1^*\cup\mathcal{S}_2 \cup\mathcal{S}_3,
$
where $\mathcal{S}_1^*=\mathcal{S}_1\backslash(\mathcal{S}_2\cup
\mathcal{S}_3)$, is disjoint. Pairs
belonging to different sets in this decomposition generate different in
law Gaussian systems.
\end{proposition}
%

\begin{pf}
We will show that Gaussian systems generated by pairs belonging to
different sets in the decomposition $\mathcal{S}=\mathcal{S}_1^*\cup
\mathcal{S}_2\cup\mathcal{S}_3$
differ by their one-dimensional distributions.
If $(\mathfrak{m},\xi)\in\mathcal{S}_2$, then $\mathfrak{m}=\alpha
\mathfrak{e}_0$ for some $\alpha>0$,
and consequently, $\mathfrak{m}_t=\mathfrak{m}*\mathfrak{n}_t=\alpha
\mathfrak{e}_0$ for every $t\in\mathbb{R}^d$.
If $(\mathfrak{m},\xi)\in\mathcal{S}_3$, then $\mathfrak{m}=\alpha
\mathfrak{e}_{\lambda}$ for some
$\alpha
>0$ and $\lambda\neq0$. Hence, in this case, $\mathfrak
{m}_t=\mathfrak{m}*\mathfrak{n}
_t=\tilde
\alpha\mathfrak{e}_{\lambda}$ for some $\tilde\alpha>0$. Finally,
if $(\mathfrak{m}
,\xi
)\in\mathcal{S}_1^*$, then $\mathfrak{m}_t$ is not a multiple of
$\mathfrak{e}_\lambda$,
$\lambda
\in\mathbb{R}$. Otherwise, Lemma~\ref{lem:conv_eq_1} would imply
that the same
is true for $\mathfrak{m}$, which contradicts the assumption
$(\mathfrak{m},\xi)\in
\mathcal{S}_1^*$.
\end{pf}

In the sequel, we concentrate on pairs belonging to the same set in the
decomposition $\mathcal{S}=\mathcal{S}_1^*\cup\mathcal{S}_2 \cup
\mathcal{S}_3$. Let us call a pair
$(\mathfrak{m},\xi)$ belonging to $\mathcal{S}_2$ or $\mathcal
{S}_3$ \textit{canonical} if
$\xi
(0)=0$. A classification of such pairs was given in Corollary~\ref
{cor:main_gen}.
\begin{proposition}\label{prop:end_class_2}
For every $(\mathfrak{m},\xi)\in\mathcal{S}_2$ there is a unique
canonical pair
$(\tilde{\mathfrak{m}}, \tilde{\xi})\in\mathcal
{S}_2$ generating the same
Gaussian system as $(\mathfrak{m},\xi)$.
\end{proposition}
\begin{pf}
To show the existence, set $\tilde\mathfrak{m}=\mathfrak{m}$ and
$\tilde\xi(t)=\xi
(t)-\xi(0)$.
Applying~(\ref{eq:wspom_br_br}) two times, we obtain that for every
$B\in\mathcal{B}(\mathbb{R}^n)$,
\begin{eqnarray*}
\mathfrak{m}_{t_1,\ldots,t_n}(B)
&=&\alpha\int_{\mathbb{R}} \mathbb{P}\bigl[\bigl(\xi(t_i)-\xi(t_1)\bigr)_{i=1}^n\in B-z\bigr]\,dz\\
&=&\alpha\int_{\mathbb{R}} \mathbb{P}\bigl[\bigl(\tilde\xi(t_i)-\tilde\xi(t_1)\bigr)_{i=1}^n\in B-z\bigr]\,dz\\
&=&\tilde\mathfrak{m}_{t_1,\ldots,t_n}(B),
\end{eqnarray*}
where $\tilde\mathfrak{m}_{t_1,\ldots,t_n}$ are the finite-dimensional
intensities of $\mathit{GS}(\tilde{\mathfrak{m}}, \tilde
{\xi})$
[cf.~(\ref{eq:mm_t_R})].
Hence, $(\mathfrak{m},\xi)$ and $(\tilde{\mathfrak{m}},
\tilde{
\xi
})$ generate equal in law Gaussian systems.

We prove the uniqueness part. Let $(\mathfrak{m},\xi)$ be a canonical
pair. Then
$\mathfrak{m}=\alpha\mathfrak{e}_0$ and $\xi(t)=W(t)+f(t)$ (see
Theorem~\ref
{theo:main_gen}). We will show that the triple $(\alpha, W, f)$ is
uniquely determined by the finite-dimensional distributions of
$\mathfrak{P}
=\mathit{GS}(\mathfrak{m},\xi)$.

First, we have $\mathfrak{m}_t=\mathfrak{m}*\mathfrak{n}_t=\alpha
\mathfrak{e}_0$ for every $t\in\mathbb{R}
^d$, and
so, $\alpha$ is uniquely determined.
Let us turn to the two-dimensional distributions of $\mathfrak{P}$.
By~(\ref{eq:mm_t_R}), we have
\[
\mathfrak{m}_{0,t}(B)=\alpha\int_{\mathbb{R}} \mathfrak{n}_{0,t}(B-z)\,dz.
\]
By Lemma~\ref{lem:proof_main}, there is a representation
\[
\mathfrak{m}_{0,t}(B)=\alpha\int_{\mathbb{R}}\mathfrak
{m}^{(0)}_{0,t}(B-z)\,dz
\]
for some measure $\mathfrak{m}^{(0)}_{0,t}$ concentrated on the line
$\{
(x_1,x_2)\in\mathbb{R}^2\dvtx x_1=0\}$ and having the Laplace transform
$\exp\{ f(t)
u + 1/2 \gamma(0,t)u^2 \}$. By Lemma~\ref{lem:proof_main_uniq}, this
shows that the two-dimensional distributions of $\mathfrak{P}$
determine $f(t)$
and $\gamma(0,t)$ uniquely. To see that $\gamma(0,t)$ determines the
law of $W$ uniquely, recall that $W(0)=0$ and hence, we may write the
covariance function of $W$ in the form
\[
r(t_1,t_2)=\tfrac12 \bigl(\gamma(0,t_1)+\gamma(0,t_2)-\gamma(0,t_1-t_2)\bigr).
\]
This completes the proof of the uniqueness part.
\end{pf}

\begin{proposition}\label{prop:end_class_3}
For every $(\mathfrak{m},\xi)\in\mathcal{S}_3$ there is a unique
canonical pair
$(\tilde
{\mathfrak{m}}, \tilde{\xi})\in\mathcal
{S}_3$ generating the same
Gaussian system as $(\mathfrak{m},\xi)$.
\end{proposition}
\begin{pf}
All necessary\vspace*{1.5pt} ingredients are contained in~\cite{kabluchkoetal}. Take
$\tilde\xi(t)=\xi(t)-\xi(0)$ and $\tilde\mathfrak{m}=\mathfrak
{m}*\delta_c$,
where $c$
is as in Theorem~\ref{theo:main_gen}. The fact that $(\tilde
{
\mathfrak{m}}, \tilde{\xi})$ and $(\mathfrak{m},\xi)$
generate equal Gaussian
systems was essentially shown in Proposition~11 of~\cite
{kabluchkoetal}. The uniqueness part follows under the additional
assumption $\lambda=1$ from Remark~24 of~\cite{kabluchkoetal}. The
general case is analogous.
\end{pf}

The next proposition gives a necessary and sufficient condition on two
pairs belonging to $\mathcal{S}_1^*$ to generate equal in law Gaussian systems.
%
\begin{proposition}\label{prop:end_class_1}
Let $(\mathfrak{m}',\xi')$ and $(\mathfrak{m}'',\xi'')$ be two
pairs, both belonging to
$\mathcal{S}_1^*$ and generating Gaussian systems $\mathfrak{P}'$ and
$\mathfrak{P}''$. Then
%
\begin{equation}\label{eq:gs_eq_gs}
\mathfrak{P}'\stackrel{d}{=}\mathfrak{P}''
\end{equation}
iff the following holds: There is a Gaussian variable $N_0$ whose
distribution on $\mathbb{R}$ is denoted by $\mathfrak{n}_0$ and which
is independent of
$\xi', \xi''$, such that
%
\begin{equation}\label{eq:gs_eq_gs_cond_1}
\mathfrak{m}'=\mathfrak{m}''*\mathfrak{n}_{0} \quad\mbox{and}\quad \{\xi
''(t), t\in\mathbb{R}^d\}\stackrel{d}{=}\{
\xi
'(t)+N_0, t\in\mathbb{R}^d\},
\end{equation}
or
%
\begin{equation}\label{eq:gs_eq_gs_cond_2}
\mathfrak{m}''=\mathfrak{m}'*\mathfrak{n}_{0}\quad \mbox{and}\quad \{\xi
'(t), t\in\mathbb{R}^d\}\stackrel{d}{=}\{
\xi
''(t)+N_0, t\in\mathbb{R}^d\}.
\end{equation}
\end{proposition}
\begin{pf}
Introduce the notation $\mu'$, $r'$, $\mu''$, $r''$, etc.\ as in
Section~\ref{subsec:notation}.
By definition of the family $\mathcal{S}_1^*$, the functions $\mu'$,
$\sigma
'^2$, $\mu''$, $\sigma''^2$ are constant. Therefore, we write, say,
$\mu
'$ instead of $\mu'(t)$. We may rewrite~(\ref{eq:r_sigma_gamma}) as follows:
%
\begin{equation}\label{eq:gamma_prime_two_prime}
\hspace*{30pt}\gamma'(t_1,t_2)=2\bigl(\sigma'^2-r'(t_1,t_2)\bigr) \quad\mbox{and}\quad \gamma
''(t_1,t_2)=2\bigl(\sigma''^2-r''(t_1,t_2)\bigr).
\end{equation}

We start by proving the ``if'' part of the proposition. Assume for
concreteness that~(\ref{eq:gs_eq_gs_cond_1}) holds. Then, by~(\ref
{eq:mm_t_R}),
\begin{eqnarray*}
\mathfrak{m}''_{t_1,\ldots,t_n}(B)
&=&\int_{\mathbb{R}}\mathbb{P}[(\xi''(t_1),\ldots,\xi''(t_n))\in
B-z]\mathfrak{m}''(dz)\\
&=&\int_{\mathbb{R}}\mathbb{P}\bigl[\bigl(\xi'(t_1)+N_0,\ldots,\xi'(t_n)+N_0\bigr)\in B-z\bigr]\mathfrak{m}''(dz)\\
&=&\int_{\mathbb{R}}\int_{\mathbb{R}}\mathbb{P}[(\xi'(t_1),\ldots,\xi'(t_n))\in B-(z+y)]\mathfrak{m}''(dz)\mathfrak{n}_0(dy).
\end{eqnarray*}
For every nonnegative function $f\dvtx \mathbb{R}^d\to\mathbb{R}$ the
following formula holds:
\[
\int_{\mathbb{R}}\int_{\mathbb{R}}f(z+y)\mathfrak
{m}''(dz)\mathfrak{n}_0(dy)=\int_{\mathbb{R}}f(x)(\mathfrak{m}
''*\mathfrak{n}_0)(dx).
\]
Hence,
\begin{eqnarray*}
\mathfrak{m}''_{t_1,\ldots,t_n}(B)
&=&\int_{\mathbb{R}}\mathbb{P}[(\xi'(t_1),\ldots,\xi'(t_n))\in
B-x](\mathfrak{m}''*\mathfrak{n}
_0)(dx)\\
&=&\int_{\mathbb{R}}\mathbb{P}[(\xi'(t_1),\ldots,\xi'(t_n))\in
B-x]\mathfrak{m}'(dx)\\
&=&\mathfrak{m}'_{t_1,\ldots,t_n}(B).
\end{eqnarray*}
This proves~(\ref{eq:gs_eq_gs}).

Now we prove the ``only if'' part of the proposition. Assume
that~(\ref{eq:gs_eq_gs}) holds. Without loss of generality we assume
that $\sigma
'^2\leq\sigma''^2$. Define
\[
\mathfrak{n}_0=\mathfrak{N}(\mu''-\mu', \sigma''^2-\sigma'^2),
\]
and let $N_0\sim\mathfrak{n}_0$ be a Gaussian variable independent of
$\xi'$ and
$\xi''$. We will show that~(\ref{eq:gs_eq_gs_cond_1}) holds.

We start by proving the first equality in~(\ref{eq:gs_eq_gs_cond_1}).
It follows from~(\ref{eq:gs_eq_gs}) that $\mathfrak{m}'_t=\mathfrak
{m}''_t$ for all
$t\in
\mathbb{R}^d$. Equivalently,
\[
\mathfrak{m}'*\mathfrak{N}(\mu',\sigma'^2)=\mathfrak
{m}''*\mathfrak{N}(\mu'',\sigma''^2).
\]
Then, by Lemma~\ref{lem:conv_eq_4}, $\mathfrak{m}'=\mathfrak
{m}''*\mathfrak{n}_{0}$. This proves
the first equality in~(\ref{eq:gs_eq_gs_cond_1}).

We claim that the second equality in~(\ref{eq:gs_eq_gs_cond_1}) follows
from the following statement: for all $t_1,t_2\in\mathbb{R}^d$,
%
\begin{equation}\label{eq:gs_eq_gs_cond_eq}
\gamma'(t_1,t_2)=\gamma''(t_1,t_2).
\end{equation}
%
To see this, set for a moment $\tilde\xi'(t)=\xi'(t)+N_0$. Then
\[
\mathbb E\tilde\xi'(t)=\mu'+(\mu''-\mu')=\mu''=\mathbb E\xi''(t).
\]
Elementary transformations using~(\ref{eq:gamma_prime_two_prime})
and~(\ref{eq:gs_eq_gs_cond_eq}) yield
\begin{eqnarray*}
\operatorname{Cov}(\tilde\xi'(t_1), \tilde\xi'(t_2))
=r'(t_1,t_2)+(\sigma''^2-\sigma'^2)
=r''(t_1,t_2)
=\operatorname{Cov}(\xi''(t_1),\xi''(t_2)).
\end{eqnarray*}

From now on, we are proving~(\ref{eq:gs_eq_gs_cond_eq}). We need to
consider two cases.
%
\begin{caseee}\label{case30psl}
Assume that $\mathfrak{m}'=\alpha'\mathfrak
{e}_{\lambda}+\beta
\mathfrak{e}
_0$ for some $\alpha'>0$, $\beta>0$, $\lambda\neq0$.
It follows from $\mathfrak{m}'_{t}=\mathfrak{m}''_{t}$ that
\[
\mathfrak{m}'*\mathfrak{N}(\mu',\sigma'^2)=\mathfrak
{m}''*\mathfrak{N}(\mu'',\sigma''^2).
\]
The left-hand side of the above equation is of the form $\alpha
\mathfrak{e}
_{\lambda}+\beta\mathfrak{e}_0$ for some $\alpha>0$. Hence, using
Lemma~\ref
{lem:conv_eq_1}, we conclude that $\mathfrak{m}''=\alpha''\mathfrak
{e}_{\lambda
}+\beta\mathfrak{e}
_0$ for some $\alpha''>0$.

Let us consider the two-dimensional distributions of $\mathfrak{P}'$.
By~(\ref{eq:mm_t_R}),
\[
\mathfrak{m}_{t_1,t_2}'(B)
=\alpha' \int_{\mathbb{R}} e^{-\lambda z}\mathfrak
{n}_{t_1,t_2}'(B-z)\,dz+\beta\int
_{\mathbb{R}
}\mathfrak{n}_{t_1,t_2}'(B-z)\,dz,\qquad B\in\mathcal{B}(\mathbb{R}^2).
\]
Applying Lemma~\ref{lem:proof_main} twice, we get two measures
$\mathfrak{m}'^{(\lambda)}_{t_1,t_2}$ and $\mathfrak{m}'^{(0)}_{t_1,t_2}$
concentrated on $\{(x_1,x_2)\in\mathbb{R}^2\dvtx  x_1=0\}$ such that the following
decomposition is valid:
\[
\mathfrak{m}'_{t_1,t_2}(B)=\alpha' \int_{\mathbb{R}}e^{-\lambda
z}\mathfrak{m}'^{(\lambda
)}_{t_1,t_2}(B-z)\,dz+
\beta\int_{\mathbb{R}}\mathfrak{m}'^{(0)}_{t_1,t_2}(B-z)\,dz,\qquad B\in
\mathcal{B}(\mathbb{R}^2).
\]
Furthermore, $\psi'_{t_1,t_2}(u)$, the Laplace transform of $\mathfrak{m}
'^{(0)}_{t_1,t_2}$,
is given by
%
\begin{equation}
\psi'_{t_1,t_2}(u)=e^{\gamma'(t_1,t_2) u^2/2}.
\end{equation}

Similar calculations can be done for $\mathfrak{m}''_{t_1,t_2}$.
By~(\ref{eq:gs_eq_gs}), we must have $\mathfrak
{m}'_{t_1,t_2}=\mathfrak{m}''_{t_1,t_2}$. By
Lemma~\ref{lem:proof_main_uniq}, this implies
\[
\mathfrak{m}'^{(0)}_{t_1,t_2}=\mathfrak{m}''^{(0)}_{t_1,t_2}.
\]
Comparing the Laplace transforms, we obtain~(\ref{eq:gs_eq_gs_cond_eq}).
\end{caseee}

\begin{caseee}
Assume that the condition of Case~\ref{case30psl} is not satisfied.
We define a point process $\tilde\mathfrak{P}'_{t_1,t_2}$ as in~(\ref
{eq:def_P_wsp110}) and (\ref{eq:def_P_wsp111}) with $\vartheta=1/2$:
we set
\begin{eqnarray*}
\tilde\mathfrak{P}'_{t_1,t_2}
= \{U_i'+ \xi'_i(t_1)/2+\xi'_i(t_2)/2, i\in\mathbb{N}\},
\end{eqnarray*}
where $\{U_i', i\in\mathbb{N}\}$ and $\xi_i'$, $i\in\mathbb{N}$,
are the starting
points and the driving processes of the Gaussian system $\mathfrak
{P}'$. Then
$\tilde\mathfrak{P}'_{t_1,t_2}$ is a Poisson\vspace*{1pt} point process on
$\mathbb{R}$ whose
intensity measure $\tilde\mathfrak{m}'_{t_1,t_2}$ is given by the formula
%
\begin{equation}\label{eq:tilde_mm_rep1}
\tilde\mathfrak{m}_{t_1,t_2}'=\mathfrak{m}'*\mathfrak{N}\bigl(\mu',
\tfrac12 \sigma'^2+\tfrac12
r'(t_1,t_2) \bigr).
\end{equation}
A simple calculation using~(\ref{eq:gamma_prime_two_prime}) shows that
\[
\mathfrak{m}_{t_1}'=\tilde\mathfrak{m}_{t_1,t_2}'*\mathfrak{N}\bigl(0,
\tfrac14 \gamma'(t_1,t_2) \bigr).
\]

Similar calculations can be done for the pair $(\mathfrak{m}'',\xi'')$.
By~(\ref{eq:gs_eq_gs}), we must have $\tilde\mathfrak
{m}'_{t_1,t_2}=\tilde
\mathfrak{m}
''_{t_1,t_2}$.
Denoting these equal measures for a moment by $\tilde\mathfrak{m}_{t_1,t_2}$,
we obtain
\begin{eqnarray*}
\tilde\mathfrak{m}_{t_1,t_2}*\mathfrak{N}\bigl(0, \tfrac14 \gamma
'(t_1,t_2) \bigr)
=\tilde\mathfrak{m}_{t_1,t_2}*\mathfrak{N}\bigl(0, \tfrac14 \gamma
''(t_1,t_2) \bigr).
\end{eqnarray*}

Now assume that~(\ref{eq:gs_eq_gs_cond_eq}) does not hold for some
$t_1,t_2\in\mathbb{R}^d$. Then Lemma~\ref{lem:conv_eq_2} implies that
$\tilde
\mathfrak{m}_{t_1,t_2}$ is of the form $\tilde\alpha\mathfrak
{e}_{\lambda}+\tilde
\beta
\mathfrak{e}_0$ for some $\tilde\alpha\geq0$, $\tilde\beta\geq0$ and
$\lambda
\neq0$. Further, Lemma~\ref{lem:conv_eq_1} applied to~(\ref
{eq:tilde_mm_rep1}) yields that $\mathfrak{m}'$ is of the form $\alpha
'\mathfrak{e}
_{\lambda}+\beta'\mathfrak{e}_0$ for some $\alpha'\geq0$, $\beta
'\geq0$ and
$\lambda\neq0$. In fact, the assumption $(\mathfrak{m}',\xi')\in
\mathcal{S}_1^*$
implies that even $\alpha'>0$, $\beta'>0$. Hence, we are in the
situation of Case~\ref{case30psl}, which is a contradiction.
\end{caseee}

The proof of Proposition~\ref{prop:end_class_1} is complete.
\end{pf}

\section{Open questions}

We have considered only particles moving on the one-dimensional real
line (although we allowed for a multidimensional time). An interesting
question is whether it is possible to obtain an analogue of
Theorem~\ref
{theo:main_gen} for particles moving in a multidimensional Euclidean space.

Another problem is to classify all stationary systems of particles
driven by independent Gaussian processes and starting at the points of
an arbitrary point process (rather than a Poisson point process).
It seems that to gain information from the stationarity of the
one-dimensional distributions of such particle systems, the results
of~\cite{liggett78} should be used instead of that of~\cite{deny59}.

\section*{Acknowledgment}

The author is grateful to Martin Schlather for several useful remarks.

%
%

\printaddresses

\end{document}